\documentclass[preprint,12pt]{elsarticle}

\usepackage{amssymb}
\usepackage{amsthm}
\usepackage{graphics}
\usepackage{epsfig}
\usepackage{graphicx}
\usepackage{color}
\usepackage{hyperref}
\usepackage{a4wide}
\usepackage{amsmath}
\usepackage{amsfonts}
\usepackage{enumerate}
\newcommand{\pf}{\noindent {\bf Proof: }}

\newtheorem{lemma}{\bf Lemma}[section]

\newtheorem{theorem}{\bf Theorem}[section]
\newtheorem{proposition}[lemma]{\bf Proposition}
\newtheorem{corollary}[lemma]{\bf Corollary}
\newtheorem{definition}{\bf Definition}[section]
\newtheorem{remark}{\bf Remark}[section]

\journal{~}

\begin{document}

\begin{frontmatter}



\title{Triameter of Graphs}



\author{Angsuman Das\corref{cor1}}
\ead{angsumandas@sxccal.edu}

\address{Department of Mathematics,\\ St. Xavier's College, Kolkata, India.\\angsumandas@sxccal.edu}
\cortext[cor1]{Corresponding author}

\begin{abstract}
In this paper, we introduce and study a new distance parameter {\it triameter} of a connected graph $G$, which is defined as $max\{d(u,v)+d(v,w)+d(u,w): u,v,w \in V\}$ and is denoted by $tr(G)$. We find various upper and lower bounds on $tr(G)$ in terms of order, girth, domination parameters etc., and characterize the graphs attaining those bounds. In the process, we provide some lower bounds of (connected, total) domination numbers of a connected graph in terms of its triameter. The lower bound on total domination number was proved earlier by Henning and Yeo. We provide a shorter proof of that. Moreover, we prove Nordhaus-Gaddum type bounds on $tr(G)$ and find $tr(G)$ for some specific family of graphs.
\end{abstract}


\end{frontmatter}


\section{Introduction}
The channel assignment problem is the problem of assigning frequencies to
the transmitters in some optimal manner and with no interferences. Keeping this problem in mind, Chartrand {\it et. al.} in \cite{radio-k} introduced the concept of {\it radio $k$-coloring} of a simple connected graph. As finding the radio $k$-chromatic number of graphs is highly non-trivial and therefore is known for very few graphs, determining good and sharp bounds is an interesting problem and has been studied by many authors \cite{panigrahi-paper},\cite{laxman-panigrahi-paper1},\cite{laxman-panigrahi-paper2},\cite{ushnish-1},\cite{ushnish-2} etc. In \cite{panigrahi-paper},\cite{laxman-panigrahi-paper1},\cite{laxman-panigrahi-paper2}, authors provides some sharp lower bounds on radio $k$-chromatic number of connected graphs in terms of a newly defined parameter called triameter of a graph (It was denoted as $M$-value of a graph in \cite{laxman-panigrahi-paper2}). Apart from this, the concept of triameter also finds application in metric polytopes \cite{metric-polytope}. Recently, in \cite{henning-total-lower-bound}, Henning and Yeo proved a graphitti conjecture on lower bound of total domination number of a connected graph in terms of its triameter. Keeping these as motivation, in this paper, we formally study triameter of connected graphs and various bounds associated with it. In fact, in the process, we provide a shorter proof of the main result in \cite{henning-total-lower-bound}.

\section{Preliminaries}
In this section, for convenience of the reader and also for later use, we recall some definitions, notations and results concerning elementary graph theory. For undefined terms and concepts the reader is referred to \cite{west-graph-book}.

By a graph $G=(V,E)$, we mean a non-empty set $V$ and a symmetric binary relation (possibly empty) $E$ on $V$.  If two vertices $u,v$ are adjacent in $G$, either we write $(u,v)\in E$ or $u \sim v$ in $G$. The distance $d_G(u,v)$ or $d(u,v)$ between two vertices $u,v\in V$ is the length of the shortest path joining $u$ and $v$ in $G$. The {\it eccentricity} of a vertex $v$ is defined as $max\{d(u,v):u\in V\}$ and is denoted by $ecc(v)$. The {\it radius}, {\it diameter} and {\it center} of a connected graph $G$ are defined as $rad(G)=min\{ecc(v):v \in V\}$, $diam(G)=max\{ecc(v):v \in V\}$ and $center(G)=\{v\in V:ecc(v)=rad(G)\}$ respectively. The Wiener index $\sigma(G)$ is defined as $\sum_{\{u,v\} \subset V}d(u,v)$. A graph $G$ is said to be vertex transitive if $\mathsf{Aut}(G)$, the automorphism group of $G$, acts transitively on $G$. The length of a cycle, if it exists, of smallest length is said to be the girth $g(G)$ of $G$. A graph $G$ is said to be Hamiltonian if there exists a cycle containing all the vertices of $G$ as a subgraph of $G$. A graph $G$ is said to be strongly regular with parameters $(n,k,\lambda,\mu)$ if it is a $k$-regular $n$-vertex graph in which any two adjacent vertices have $\lambda$ common neighbours and any two non-adjacent vertices have $\mu$ common neighbours. A graph is said to be a bistar if it is obtained by joining the root vertices of two stars $K_{1,n_1}$ and $K_{1,n_2}$. We denote this graph by $K^{n_1}_{n_2}$ and it is a graph on $n_1+n_2+2$ vertices.

\section{Triameter of a Graph and its Bounds}
In what follows, even if not mentioned, $G$ denotes a finite simple connected undirected graph with at least $3$ vertices. We start by defining triameter of a connected graph.
{\definition Let $G=(V,E)$ be a connected graph on $n\geq 3$ vertices. The triameter of $G$ is defined as $\max\{d(u,v)+d(v,w)+d(u,w): u,v,w \in V\}$ and is denoted by $tr(G)$.}

From the definition, it follows that $tr(G)$ is always greater than or equal to $3$. However, triameter of a graph on $n$ vertices can be as large as $2n-2$, as evident from the following results proved in \cite{panigrahi-paper}: $tr(P_n)=2(n-1)$ and $tr(C_n)=n$. 

If $G$ and $H$ be two connected graphs on same vertex set with $E(H)\subseteq E(G)$, then by definition of triameter, we have $tr(G)\leq tr(H)$. For any three vertices $u,v,w$, let us denote by $d(u,v,w)$, the sum $d(u,v)+d(v,w)+d(u,w)$. Now, we investigate other bounds on $tr(G)$. 
{\theorem \label{diameter-bound} For any connected graph $G$, $2\cdot diam(G)\leq tr(G) \leq 3\cdot diam(G)$ and the bounds are tight.}\\
\pf The upper bound follows from the definition of diameter and triameter of a connected graph. For the lower bound, let $d(u,v)=diam(G)$. Choose $w \in V \setminus \{u,v\}$. Then $d(u,v)\leq d(v,w)+d(w,u)\Rightarrow 2\cdot diam(G)=2d(u,v)\leq d(u,v)+d(v,w)+d(w,u)\leq tr(G)$.

The tightness of the bounds follows from the following examples: For $n\geq 3$, $tr(P_n)=2\cdot diam(P_n)$. For Petersen graph $P$, $tr(P)=3\cdot diam(P)$. \qed

{\corollary Let $G$ be a connected graph on $n$ vertices such that $\delta(G)\geq \frac{n}{2}$. Then $tr(G)\leq 6$.}\\
\pf It follows from Theorem \ref{diameter-bound} and the fact that $\delta(G)\geq \frac{n}{2}$ implies $diam(G)\leq 2$.\qed


{\corollary \label{radius-bound} For any connected graph $G$, $2\cdot rad(G) \leq tr(G) \leq 6\cdot rad(G)$ and the bounds are tight.}\\
\pf As for any connected graph $G$, $rad(G)\leq diam(G)\leq 2\cdot rad(G)$, we have $2\cdot rad(G) \leq tr(G) \leq 6\cdot rad(G)$. For the tightness of lower bound, take $G=C_{2n}$ where $tr(G)=2n=2\cdot rad(G)$ and for upper bound, take $G=K_{1,3}$ where $tr(G)=6$ and $rad(G)=1$.\qed

\begin{figure}[h]
	\centering
	\begin{center}
		\includegraphics[scale=.35]{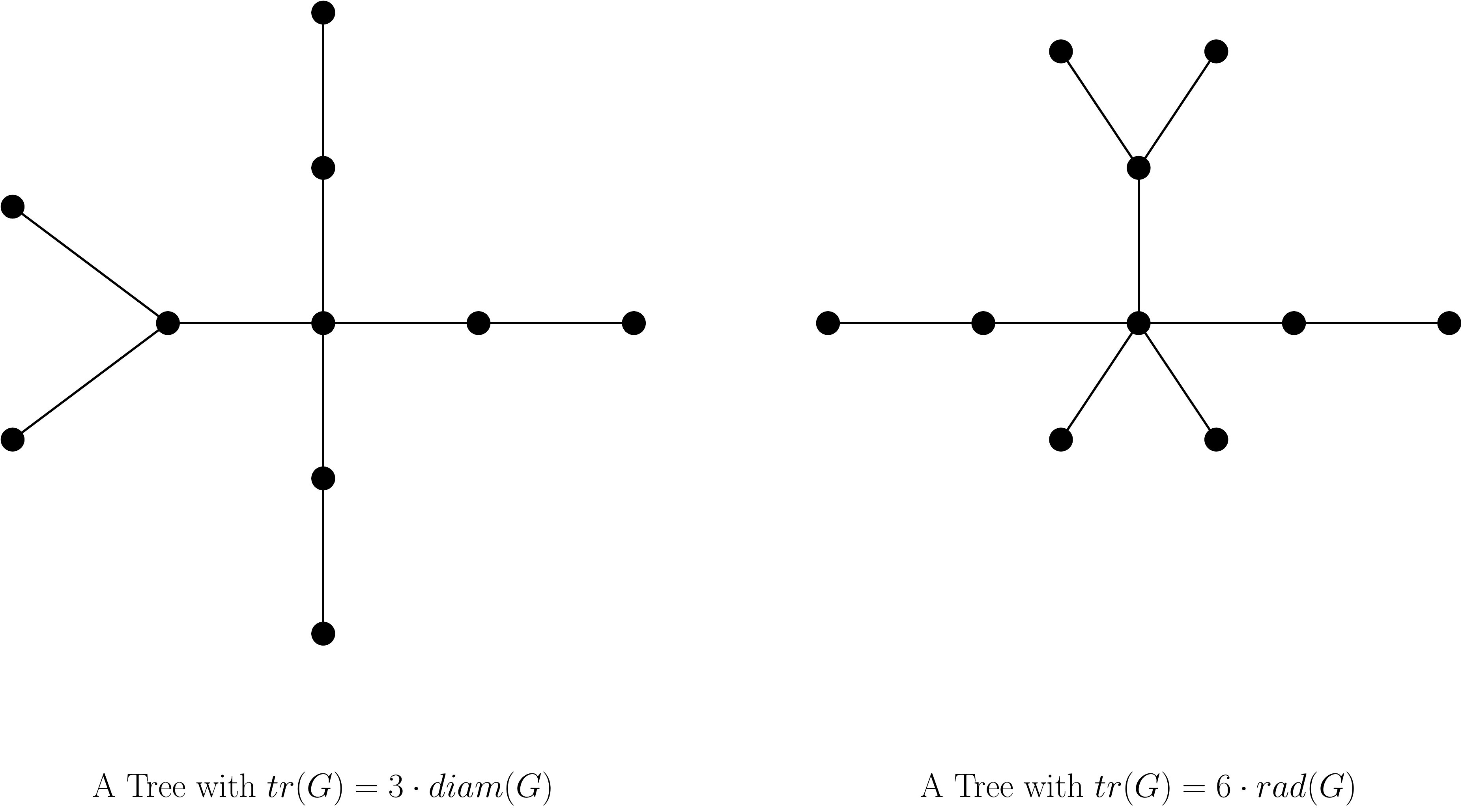}
		\caption{Trees achieving the Upper Bounds}
		\label{tree-figure}
	\end{center}
\end{figure}

{\remark Some other examples demonstrating the tightness of the upper bounds are shown in Figure \ref{tree-figure}. The bound in Corollary \ref{radius-bound} can be substantially tightened in case of vertex transitive graphs. See Theorem \ref{vt-bound}.}

{\corollary \label{tree-radius-bound} For any tree $T$, $4\cdot rad(T)-2 \leq tr(T) \leq 6\cdot rad(T)$ and the bounds are tight.}\\
\pf We first recall a result on tree: A tree $T$ has either $|center(T)|=1$ or $|center(T)|=2$, and $diam(T)=2\cdot rad(T)$ or $2\cdot rad(T)-1$ according as $|center(T)|=1$ or $|center(T)|=2$. Hence the corollary follows from Theorem \ref{diameter-bound}. Tightness of upper bound and lower bound follows respectively from $K_{1,3}$ and $P_4$. \qed

It is known that in a connected graph $G$ with cycle, $g(G)\leq 2\cdot diam(G)+1$. Thus it trivially follows from Theorem \ref{diameter-bound} that $g(G)\leq tr(G)+1$. In the next theorem, we prove a stronger inequality involving girth and triameter.

\subsection{Upper Bounds}

{\theorem \label{order-bound} For any connected graph $G$ with $n\geq 3$ vertices, $tr(G)\leq 2n-2$ and the bound is tight.}\\
\pf It suffices to prove the bound for trees, as for any connected graph $G$ and any spanning tree $T$ of $G$, $tr(G)\leq tr(T)$ holds. We prove the result by induction on $n$. Clearly, in the basis step, $n=3$ and there exists only one tree on $3$ vertices, i.e., $P_3$ and the result holds for $P_3$. Let the result be true for all trees with order $n-1$ and $T$ be a tree of order $n$. 

Let $x$ be an arbitrary pendant vertex $T$ and $T'$ be the tree obtained by deleting $x$ from $T$. Thus $T'$ is a tree with $n-1$ vertices and by induction hypothesis, $tr(T')\leq 2(n-1)-2=2n-4$. Let $u,v,w$ be three distinct arbitrary vertices in $T$. 

{\bf Case 1:} If none of them coincides with $x$, then $$d_T(u,v)+d_T(v,w)+d_T(w,u)=d_{T'}(u,v)+d_{T'}(v,w)+d_{T'}(w,u)$$
$$~~~~~~~~~~~~~~~~~~~~~~~~~~~~~~~~\leq tr(T')\leq 2n-4 \leq 2n-2.$$ 

{\bf Case 2:} If one of $u,v,w$ coincides with $x$, say $x=u$, then $$d_T(u,v)=1+d_{T'}(y,v)$$ $$d_T(u,w)=1+d_{T'}(y,w)$$ $$d_T(v,w)=d_{T'}(v,w),$$ where is the support vertex for $x$ in $T$. Adding the above three equations, we get $$d_T(u,v)+d_T(v,w)+d_T(w,u)=2+d_{T'}(y,v)+d_{T'}(y,w)+d_{T'}(v,w)$$
$$~~~~~~~~~~~~~~~~~~~~\leq 2 + tr(T')\leq 2n-2.$$ Combining both the cases, we get $tr(T)\leq 2n-2$ and hence the result follows by induction.

The tightness of the bound is achieved for paths, as $tr(P_n)=2n-2$. \qed

{\theorem \label{order-equality} For any connected graph $G$ with $n\geq 3$ vertices, $tr(G)=2n-2$ if and only if $G$ is a tree with $2$ or $3$ leaves.}\\
\pf If $G$ is a tree on $n$ vertices with $2$ leaves, then $G=P_n$ and $tr(G)=2n-2$. Let $G$ be a tree on $n$ vertices with $3$ leaves $u^*$, $v^*$ and $w^*$. Clearly $G$ is obtained by subdividing the edges of $K_{1,3}$ and $G$ has a unique vertex, say $x$, of degree $3$. Note that $d(u,v,w)$ is maximized for $(u,v,w)=(u^*,v^*,w^*)$. Let $d(u^*,v^*)=k_1+k_2$ where $d(u^*,x)=k_1$ and $d(x,v^*)=k_2$ and let $d(x,w^*)=k_3$. Then by counting the number of vertices in $G$, we have $n=k_1+k_2+k_3+1$. Thus $tr(G)=d(u^*,v^*,w^*)=(k_1+k_2)+(k_1+k_3)+(k_2+k_3)=2(k_1+k_2+k_3)=2(n-1)$. Thus for trees with $2$ or $3$ leaves, $tr(G)=2n-2$ holds.

Conversely, let $G$ be a connected graph on $n$ vertices with $tr(G)=2n-2$. First we show that $G$ can not be a tree with more than $3$ leaves. 

Let $G$ be a tree on $n$ vertices with $l>3$ leaves. Then $d(u,v,w)$ is maximized for $(u,v,w)=(u^*,v^*,w^*)$ for $3$ leaves $u^*,v^*,w^*$ suitably chosen from $l$ leaves, i.e., $tr(G)=d_G(u^*,v^*,w^*)$ . Let $P_1,P_2,P_3$ be the shortest paths joining the pairs $(u^*,v^*), (v^*,w^*)$ and $(w^*,u^*)$ in $G$. Since $G$ is a tree and $u^*,v^*,w^*$ are leaves, $T=P_1\cup P_2\cup P_3$ is a tree. Here union of paths denote the subgraph induced by the vertices in $T$. Now the number of vertices in $T$ (say $k$) is less than $n$, as other $l-3$ leaves of $G$ are not in $T$. Also, $T$ is a tree with exactly $3$ leaves. Hence by the previous argument for the case of $3$ leaves, we get $tr(G)=d_G(u^*,v^*,w^*)=d_T(u^*,v^*,w^*)=tr(T)=2k-2<2n-2$. Thus, for trees with more than $3$ leaves, $tr(G)<2n-2$. Thus, $G$ can not be a tree with more than $3$ leaves. 

Next we show that $G$ can not be a connected graph which is not a tree.
Let, if possible, $G$ be a connected graph with cycles and $tr(G)=2n-2$. If it has a spanning tree $T$ with more than $3$ leaves, then $tr(G)\leq tr(T)<2n-2$, a contradiction. Thus all the spanning trees of $G$ must have $2$ or $3$ leaves. Let $T$ be a spanning tree of $G$ with $2$ or $3$ leaves. If $T$ has $2$ leaves, then $T=P_n$. As $G$ contains cycle and the vertex set for $G$ and $T$ are same, $tr(G)<tr(T)=tr(P_n)=2n-2$, a contradiction. Thus let us assume that $T$ has $3$ leaves, say $u^*,v^*,w^*$. Also, let $tr(G)$ be attained by the vertices $\hat{u},\hat{v},\hat{w}$ of $G$. Note that $\hat{u},\hat{v},\hat{w}$ may not be same as $u^*,v^*,w^*$. Now two cases may arise.

{\bf Case 1:} $\{\hat{u},\hat{v},\hat{w}\}=\{u^*,v^*,w^*\}$. Then we have $$tr(G)=d_G(\hat{u},\hat{v},\hat{w})= d_G(u^*,v^*,w^*)< d_T(u^*,v^*,w^*)=tr(T)=2n-2,$$  a contradiction. In this case, the strict inequality holds as $G$ contains cycles.

{\bf Case 2:} $\{\hat{u},\hat{v},\hat{w}\}\neq \{u^*,v^*,w^*\}$. Then we have $$tr(G)=d_G(\hat{u},\hat{v},\hat{w})\leq d_T(\hat{u},\hat{v},\hat{w})< d_T(u^*,v^*,w^*)=tr(T)=2n-2,$$ a contradiction. In this case, the strict inequality holds as $d_T$ has a unique maximum at $(u^*,v^*,w^*)$.\qed


{\theorem \label{2n-2l+4-bound} Let $T$ be a tree on $n\geq 3$ vertices and $l\geq 4$ leaves. Then $tr(T)\leq 2n-2l+4$.}\\
\pf Let $tr(T)=d(u^*,v^*,w^*)$ for three leaves $u^*,v^*,w^*$ of $T$. Let $T'$ be the tree on $n-(l-3)$ vertices obtained by deleting the remaining $l-3$ leaves from $T$. Thus $tr(T)=tr(T')\leq 2(n-l+3)-2=2n-2l+4$, by Theorem \ref{order-bound}.  \qed
 
{\corollary \label{l=4} Let $T$ be a tree on $n\geq 3$ vertices such that $tr(T)=2n-4$, then $T$ has exactly $4$ leaves.}\\
\pf From Theorem \ref{2n-2l+4-bound}, we get $2n-4=tr(T)\leq 2n-2l+4$, i.e., $l\leq 4$. If $l=2$ or $3$, then $tr(T)=2n-2\neq 2n-4$. Thus $l=4$.\qed

\begin{figure}[h]
	\centering
	\begin{center}
		\includegraphics[scale=.3]{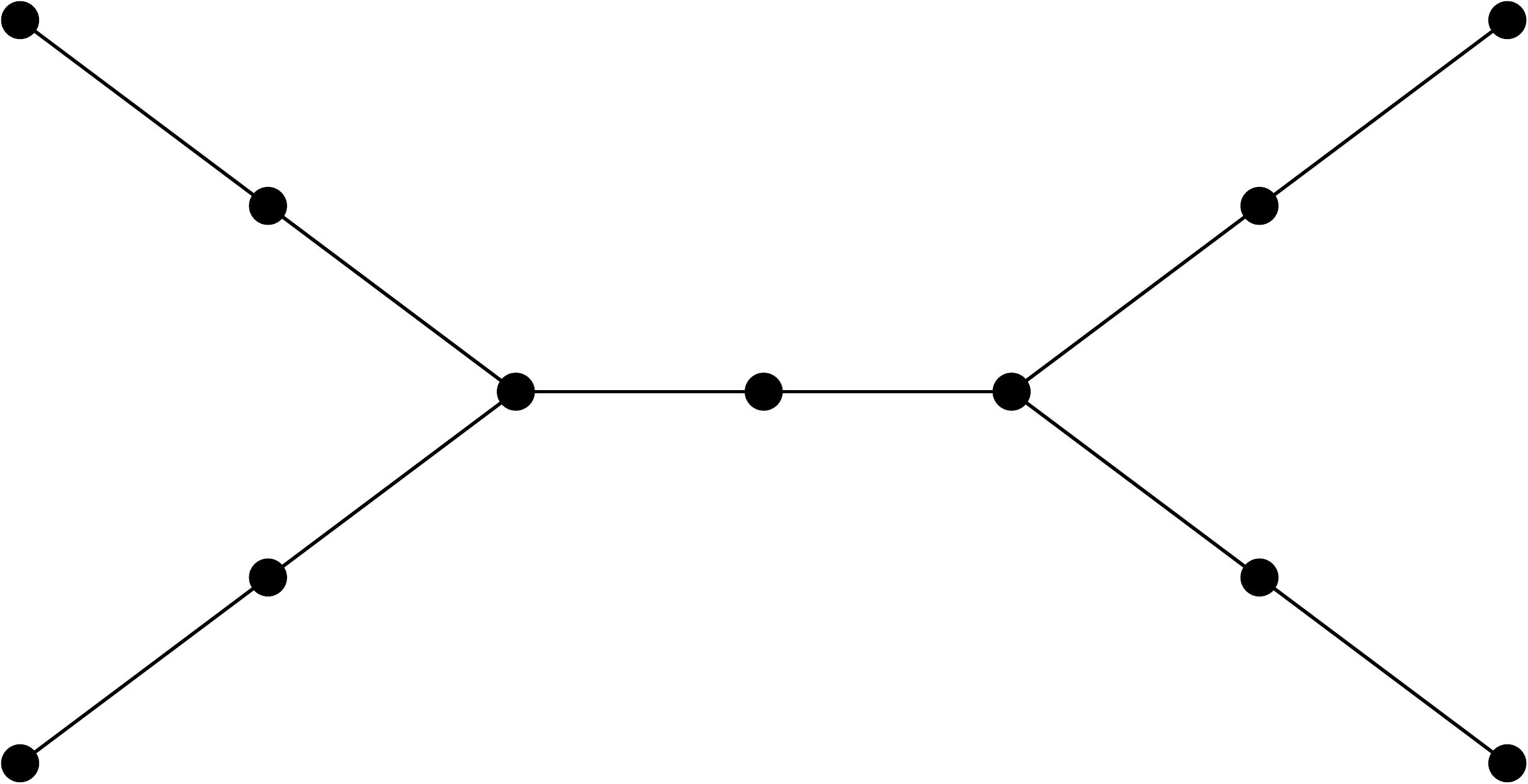}
		\caption{Example of a tree $T$ with $l(T)=4,n(T)=11$ and $tr(T)=16<18=2\cdot 11-4$ }
		\label{l=4-figure}	
	\end{center}
\end{figure}

It is to be noted that the converse of the above corollary is not true. See Figure \ref{l=4-figure}.
{\corollary \label{conn-dom-bound} Let $G$ be a connected graph on $n$ vertices with connected domination number $\gamma_c$. Then $tr(G)\leq 2\gamma_c+4$.}\\
\pf Let $T$ be a spanning tree of $G$ with maximum number of leaves $l$. Then $l+\gamma_c=n$. Now, if $l\geq 4$, $tr(G)\leq tr(T)\leq 2(n-l)+4=2\gamma_c+4$. If $l=2$ or $3$, by Theorem  \ref{order-equality}, $tr(G)=2n-2$ and $\gamma_c=n-2$ or $n-3$. In this case also, $tr(G)\leq 2\gamma_c+4$ holds.\qed

{\corollary \label{dom-bound} Let $G$ be a connected graph with domination number $\gamma(G)$. Then $tr(G)\leq 6\gamma(G)$ and the bound is tight.}\\
\pf It follows from the fact that $tr(G)\leq 2\gamma_c(G)+4$ and $\gamma_c(G)\leq 3\gamma(G)-2$ (See \cite{gamma-and-gamma_c-connection}). The bound is achieved by $K_{1,n}$. \qed

{\corollary \label{total-dom-bound} Let $G=(V,E)$ be a connected graph  with total domination number $\gamma_t(G)$. Then $tr(G)\leq 4\gamma_t(G)$.}\\
\pf In \cite{gamma_c-and-gamma_t-relation}, it was shown that $\gamma_c(G)\leq 2\gamma_t(G)-2$. Thus from Corollary \ref{conn-dom-bound}, we get $tr(G)\leq 2\gamma_c+4 \leq 2(2\gamma_t(G)-2)+4\leq 4\gamma_t(G)$.\qed

{\remark Corollary \ref{total-dom-bound} was also proved in \cite{henning-total-lower-bound}. However, here we provide a shorter proof of $tr(G)\leq 4\gamma_t(G)$ using Theorems \ref{order-bound} and \ref{2n-2l+4-bound} and Corollaries \ref{conn-dom-bound} and \ref{total-dom-bound}.}

In the next proposition, we show that the upper bound proved in Theorem \ref{order-bound} can be substantially tightened if the vertex connectivity $\kappa$ of $G$ increases.

{\proposition Let $G$ be a graph on $n$ vertices with vertex connectivity $\kappa$. Then $tr(G)\leq \dfrac{3(n-2)}{\kappa}+3$.}\\
\pf The proof follows from the result that $n \geq \kappa (diam(G)-1)+2$ (See Pg 174, Sum no. 4.2.22, \cite{west-graph-book}) and $tr(G)\leq 3\cdot diam(G)$.\qed

{\theorem \label{chromatic-corollary} For a connected graph $G$, other than odd cycle and complete graph, on $n$ vertices with maximum degree $\Delta(G)$ and chromatic number $\chi(G)$,  $tr(G)+\chi(G)\leq tr(G)+\Delta(G)\leq 2n+1$, with equality holding only if $G$ is a tree with $3$ leaves. However, for odd cycles and complete graphs, $tr(G)+\Delta(G)< tr(G)+\chi(G)\leq 2n$.}\\
\pf We first observe that the result holds for odd cycles and complete graphs, i.e., for $G=C_n$ with odd $n\geq 5$, $tr(G)=n, \chi(G)=3, \Delta(G)=2$ and for $G=K_n$, $tr(G)=3, \chi(G)=n, \Delta(G)=n-1$. Thus, we assume that $G$ is neither an odd cycle nor a complete graph. Let $T$ be a spanning tree of $G$ with maximum degree $\Delta(T)=\Delta(G)$. Also, the number of leaves $l(T)$ of $T$ satisfies $\Delta(T)\leq l(T)$. Therefore, by Brooks' Theorem, we have 

\begin{equation}\label{chromatic-equation}
tr(G)+\chi(G)\leq tr(G)+\Delta(G)\leq tr(T)+\Delta(T)\leq tr(T)+l(T).
\end{equation}

Now, if $l(T)\geq 4$, then by Theorem \ref{2n-2l+4-bound}, $tr(T)+l(T)\leq 2n-l+4\leq 2n$. If $l(T)=2$, then by Theorem \ref{order-equality}, $tr(T)+l(T)= 2n-2+2=2n$. If $l(T)=2$, then by Theorem \ref{order-equality}, $tr(T)+l(T)= 2n-2+3=2n+1$. Combining all the cases for $l(T)$ in Equation \ref{chromatic-equation}, we get $tr(G)+\chi(G)\leq tr(G)+\Delta(G)\leq 2n+1$.

Observe that the number of leaves $l(T)$ of $T$ satisfies $\Delta(T)\leq l(T)$ with equality holding only if $T$ is obtained by subdividing the edges of the star $K_{1,\Delta}$. Thus, if $T$ is not obtained by subdividing the star $K_{1,\Delta}$, then $tr(T)+\Delta(T)< tr(T)+l(T)\leq 2n+1$, i.e., $tr(G)+\chi(G)\leq tr(G)+\Delta(G)\leq 2n$.

Now, let us assume that $T$ is obtained by subdividing the star $K_{1,\Delta}$, i.e., $\Delta(T)=l(T)$. 

If $\Delta(T)= l(T)=2$ or $\geq 4$, $tr(T)+\Delta(T)\leq 2n$.

If $\Delta(T)= l(T)=3$, $tr(T)+\Delta(T)=2n+1$, i.e., the bound is tight for trees with $l=3$. Now, let $G$ be a connected graph which is not a tree itself, but for each of its spanning tree $T$, $\Delta(T)= l(T)=3$ holds. Since we can always choose a spanning tree $T$ of $G$ such that $\Delta(T)=\Delta(G)$, we have $\Delta(T)=\Delta(G)=l(T)=3$. Also as $G$ is not a tree, $tr(G)\leq 2n-3$. Thus $tr(G)+\chi(G)\leq tr(G)+\Delta(G)\leq 2n-3+3=2n$. Thus, equality holds only if $G$ is a tree with $3$ leaves. \qed

\subsection{Lower Bounds}

{\theorem \label{girth-bound} If $G$ be a connected graph with cycles, then $g(G)\leq tr(G)$.}\\
\pf Let $C$ be a cycle of length $g(G)=g$. Since $C$ is a smallest cycle in $G$, there exists two vertices $u$ and $v$ on $C$ such that $d(u,v)=\lfloor g/2\rfloor$. Choose $w$ on $C$ such that $w\sim v$ and $d(u,w)=\lfloor (g-1)/2\rfloor$. Again, existence of such $w$ is guaranteed as $C$ is a smallest cycle in $G$. Now, $d(u,v)+d(u,w)+d(v,w)=\lfloor g/2\rfloor + \lfloor (g-1)/2\rfloor +1 = (g-1)+1=g$ and hence the bound follows. \qed

{\theorem \label{girth=triameter-characterization} In a connected graph $G$, $g(G)=tr(G)$ holds if and only if $G$ is a complete graph or a cycle.}\\
\pf It is clear that if $G$ is a cycle, then $tr(G)=g(G)=$ length of the cycle and if $G$ is a complete graph $K_n$ with $n\geq 3$, then $tr(G)=g(G)=3$. Conversely, let $tr(G)=g(G)$ holds for a graph $G$. If $tr(G)=g(G)=3$, then $d(u,v)=1$ for all vertices $u,v$ in $G$, i.e., $G$ is a complete graph $K_n$. Also, as $g(G)=3$, we have $n\geq 3$. Thus let $tr(G)=g(G)>3$ and $C$ be a cycle of length $g=g(G)$ in $G$. Since $C$ is a smallest cycle, $C$ is a chordless induced cycle in $G$. If $G=C$, then the proof is over. If not, let $v$ be a vertex in $G$, but not in $C$, which is adjacent to some vertex $u$ in $C$, i.e., $d(u,v)=1$. 

{\bf Case 1:} $g$ is odd, say $g=2k+1>3$, i.e., $k>1$. Then there exist two vertices $x$ and $y$ in $C$ such that $d(u,x)=k=d(u,y)$ and $d(x,y)=1$. Since the girth is $2k+1$, $d(v,x)$ and $d(v,y)$ are greater or equal to $k$, otherwise we get a cycle of length less than $g$. If any one of them is greater than $k$, say $d(v,x)>k$, then we get $d(u,v)+d(u,x)+d(v,x)>1+k+k=2k+1=g$, i.e., $tr(G)>g$, a contradiction. Thus, let both $d(v,x)=d(v,y)=k$. 
Since $C$ is a cycle, there are two vertices in $C$ which are adjacent to $x$, one being $y$. Let the other vertex in $C$ which is adjacent to $x$ be $z$. Thus $d(y,z)\leq 2$ via a path through $x$. However, as $C$ is chordless, $d(y,z)\neq 1$. Thus $d(y,z)=2$. Also $d(u,z)= k-1$ as $d(u,x)=k$. Now if $d(v,z)<k$, then we get a cycle of length less than $2k+1$ passing through $u$, $v$ and $z$. Thus $d(v,z)=k$ via a path through $u$. Hence, $$tr(G)\geq d(v,z)+d(v,y)+d(y,z)=k+k+2>2k+1=g(G), \mbox{ a contradiction.}$$

\begin{figure}[h]
	\centering
	\begin{center}
		\includegraphics[scale=.3]{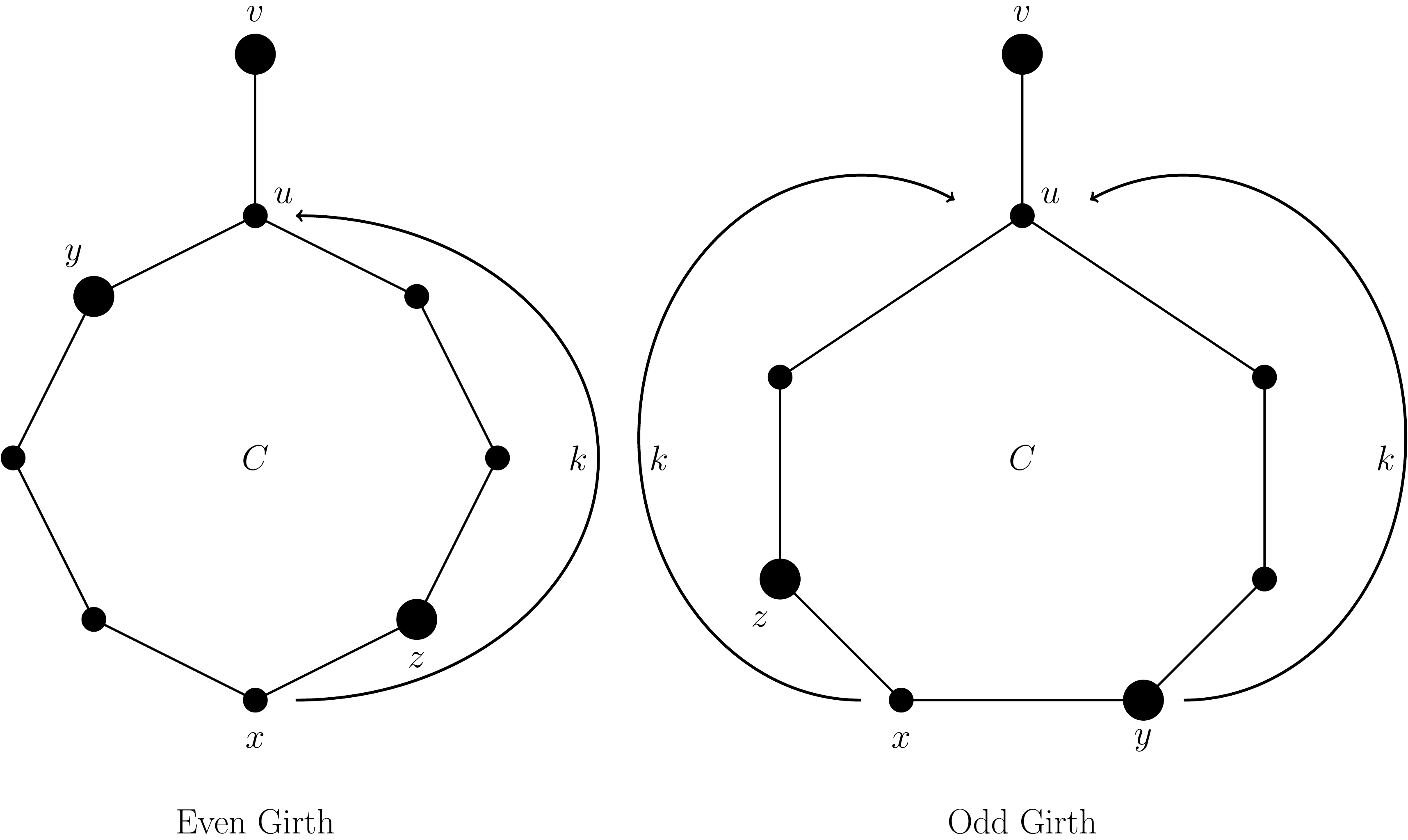}
		\caption{Schematic Diagram of the proof of Theorem \ref{girth=triameter-characterization}}
		\label{girth-figure}
	\end{center}
\end{figure}

{\bf Case 2:} $g$ is even, say $g=2k>3$, i.e., $k>1$. Then there exist a unique vertex $x$ in $C$ such that $d(u,x)=k$. Let $y$ be a vertex in $C$ adjacent to $u$. Since $k>1$, $y \neq x$. Similarly let $z$ be the unique vertex in $C$ such that $d(y,z)=k$. Note that as $C$ is a smallest cycle in $G$, $d(x,z)=1$ and $d(u,z)=k-1$. Again, $d(v,z)\geq k-1$, because if $d(v,z)< k-1$, we get a cycle of length less than $k$ through $v$, $u$ and $z$ in $G$, a contradiction. Also, $d(y,v)=2$. Hence 
$$tr(G)\geq d(v,z)+d(v,y)+d(y,z)\geq (k-1)+2+k=2k+1>g(G), \mbox{ a contradiction.}$$

Thus, combining both the cases,  there does not exist any vertex $v$ in $G$ which is not in $C$. Moreoer, as $C$ is an induced chordless cycle in $G$, we have $G=C$, i.e., $G$ is a cycle.
\qed

{\theorem \label{(n-1/l-1)-bound} Let $T$ be a tree on $n$ vertices with $l\geq 3$ leaves. Then $tr(T)\geq \left\lceil \dfrac{4(n-1)}{(l-1)} \right\rceil$ and the bound is tight. }\\
\pf For $l=3$,  its an equality. So we assume that $l>3$. Let $tr(T)=d(u,v,w)$ for three leaves $u,v,w$ in $T$. Let $P_1,P_2,P_3$ be the unique shortest path joining $u-v$,$v-w$ and $w-u$ respectively. Let $T'=\langle P_1\cup P_2 \cup P_3 \rangle$ be the sub-tree of $T$ induced by the union of $P_1,P_2$ and $P_3$. Note that $T'$ is a tree of with three leaves $u,v,w$ and $tr(T')=tr(T)$. As $T'$ is a tree with $3$ leaves, it is obtained by subdividing edges of $K_{1,3}$. Let $y$ be the root vertex in $T'$. Let $d(u,y)=k_1,d(v,y)=k_2$ and $d(w,y)=k_3$. Then $tr(T)=tr(T')=2(k_1 +k_2 +k_3)$.

\begin{figure}[h]
	\begin{center}
		\includegraphics[scale=.45]{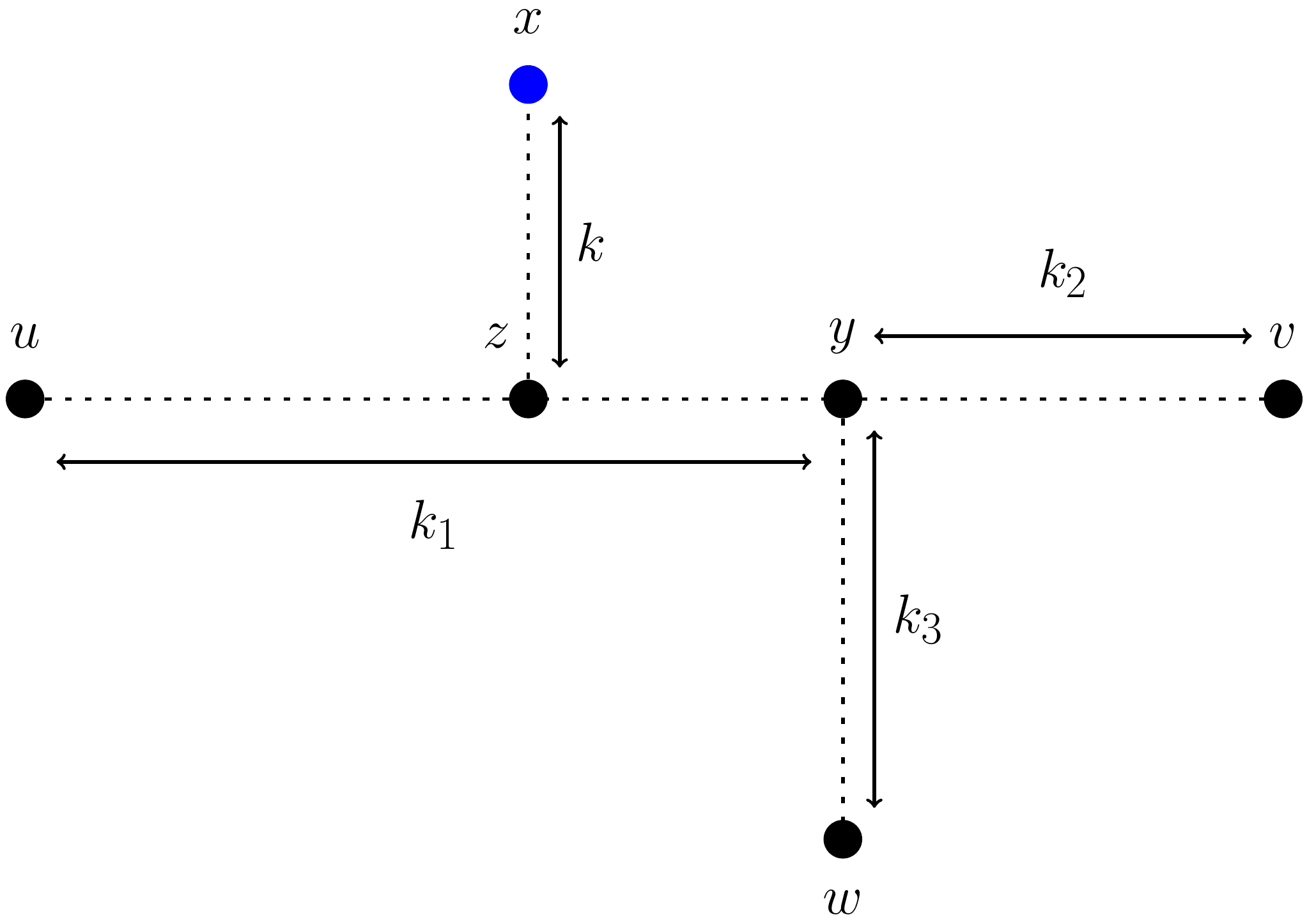}
		\caption{Schematic Diagram for proof of Theorem \ref{(n-1/l-1)-bound}}
		\label{lower-bound-figure}
	\end{center}
\end{figure}

Since, $l>3$, let $x$ be another leaf in $T$ apart from $u,v,w$ and $d(x,T')=k$, i.e., there exists $z \in T'$ such that $d(x,z)=k$ and $d(x,t)>k$ for all $t \in T'\setminus\{z\}$. Without loss of generality, let $z$ lie on the path joining $u$ and $y$. See Figure \ref{lower-bound-figure}. Here the black vertices denote the vertices of $T'$ and the blue vertex is $x$.

Claim 1: $d(u,z)\geq d(x,z)=k$.\\
If possible, let $d(u,z)< d(x,z)$, then $$d(u,v)=d(u,z)+d(z,v)<d(x,z)+d(z,v)=d(x,v) \mbox{ and}$$ $$d(u,w)=d(u,z)+d(y,z)+d(y,w)<d(x,z)+d(y,z)+d(y,w)=d(x,w).$$ $$\mbox{Thus }tr(T)=d(u,v,w)=d(u,v)+d(u,w)+d(v,w)~~~~~~~~~~~~~~~~~~~~~~~~~~~~~~~~~~~$$ $$~~~~~~~<d(x,v)+d(x,w)+d(v,w)=d(x,v,w), \mbox{ a contradiction.}$$ 

Claim 2: Either $d(v,z)\geq d(x,z)$ or $d(w,z)\geq d(x,z)$.\\
If possible, let $d(v,z)<d(x,z)$ or $d(w,z)< d(x,z)$. Without loss of generality, let $k_2\leq k_3$. Then
$$d(u,v,w)=d(u,v)+d(w,u)+d(v,w)=\left(d(u,z)+d(z,v)\right)+(k_2+k_3)+d(w,u)$$
$$<d(u,y)+d(x,z)+(k_2+k_3)+d(w,u)\mbox{ [as, }d(u,z)\leq d(u,y);d(v,z)<d(x,z)\mbox{]}$$
$$=(d(u,z)+d(y,z))+d(x,z)+(k_2+k_3)+d(w,u)~~~~~~~~~~~~~~~~~~~~~~~~~~~~~~~~~$$
$$=(d(u,z)+d(x,z))+(d(y,z)+k_3)+k_2+d(w,u)~~~~~~~~~~~~~~~~~~~~~~~~~~~~~~~~~$$
$$=d(u,x)+d(w,z)+k_2+d(w,u)~~~~~~~~~~~~~~~~~~~~~~~~~~~~~~~~~~~~~~~~~~~~~~~~~~~~~~~$$
$$<d(u,x)+d(x,z)+k_2+d(w,u)~~~~~~~~~~~~~~~~~~~~~\mbox{ [as, }d(w,z)<d(x,z)\mbox{]}~~~~~$$
$$\leq d(u,x)+(d(x,y)+k_3)+d(w,u)~~~~~~~~~~~~~\mbox{ [as, }d(x,z)\leq d(x,y);k_2\leq k_3\mbox{]}$$
$$=d(u,x)+d(x,w)+d(w,u)=d(u,x,w), \mbox{ a contradiction.}~~~~~~~~~~~~~~~~~~~~~$$

As $d(x,z)=k$, from the above two claims, we have $d(u,z)\geq k$ and either $d(v,z)$ or $d(w,z)\geq k$. Thus adding them, we get $d(u,z)+d(v,z)\geq 2k$ or $d(u,z)+d(w,z)\geq 2k$, i.e., $d(u,y)+d(v,y)=k_1+k_2\geq 2k$ or $d(u,y)+d(w,y)=k_1+k_3\geq 2k$. In any case, $2k\leq k_1+k_2+k_3$, i.e., 
\begin{equation}\label{equation-1}
k\leq \dfrac{k_1+k_2+k_3}{2}\leq \dfrac{tr(T')}{4}=\dfrac{tr(T)}{4}.
\end{equation} Let $n'$ be the number of vertices in $T'$. Then $$n'=(k_1+1)+(k_2+1)+(k_3+1)-2=k_1+k_2+k_3+1=\dfrac{tr(T)}{2}+1.$$

From Equation \ref{equation-1}, we note that while deleting vertices from $T$ to get $T'$, we have deleted at most $\dfrac{tr(T)}{4}(l-3)$ vertices, i.e., $$\dfrac{tr(T)}{2}+1=n'\geq n-\dfrac{tr(T)}{4}(l-3)$$
$$\mbox{i.e., }2tr(T)+4\geq 4n-(l-3)tr(T)$$
$$\mbox{i.e., }(l-1)tr(T)\geq 4(n-1)\Rightarrow tr(T)\geq \left\lceil \dfrac{4(n-1)}{(l-1)} \right\rceil .$$ 
The lower bound is achieved by any tree with $3$ leaves. \qed

{\theorem \label{wiener-bound} Let $G=(V,E)$ be a connected graph on $n$ vertices with Wiener index $\sigma$. Then $tr(G)\geq \dfrac{6\sigma}{n(n-1)}$ and the bound is tight.}\\
\pf Observe that for any pair of vertices $u,v \in V$, $d(u,v)$ appears $\binom{n-2}{1}$ times in the sum $\sum_{\{u,v,w \}\subset V}d(u,v,w)$. Thus, we get
$$\binom{n}{3}\cdot tr(G)\geq \sum_{\{u,v,w\}\subset V}d(u,v,w)=\binom{n-2}{1}\sum_{\{u,v\}\subset V}d(u,v)=(n-2)\sigma,$$ and hence the theorem follows. The tightness of the bound follows by taking $G=C_4$, the cycle on $4$ vertices for which $\sigma=8,tr(G)=4. $\qed

\section{Triameter of Some Graph Families}
In this ection, we find the triameter of some important families of graphs. We start by recalling a result from \cite{panigrahi-paper}.
{\proposition \cite{panigrahi-paper} For any two connected graphs $G$ and $H$, $tr(G \square H)=tr(G)+tr(H)$.}

{\corollary Let $G$ be a $m \times n$ rectangular grid graph. Then $tr(G)=2(m+n-2)$.}\\
\pf Since $G$ is a $m \times n$ rectangular grid graph, $G \cong P_m \square P_n$. Thus $tr(G)=tr(P_m)+tr(P_n)=(2m-2)+(2n-2)=2(m+n-2)$.\qed

{\theorem \label{bipartite-even-triameter} Let $G$ be a connected bipartite graph. Then $tr(G)$ is even. }\\
\pf Let $V(G)=V_1\cup V_2$ be the bipartition and $u,v,w$ be $3$ vertices in $V(G)$ such that $tr(G)=d(u,v,w)$. If $u,v,w \in V_i$ for same $i$, then each of $d(u,v),d(v,w)$ and $d(w,u)$ are even and hence $tr(G)$ is even. Thus, without loss of generality, let $u,v \in V_1$ and $w\in V_2$. Then $d(u,w)$ and $d(v,w)$ are odd and $d(u,v)$ is even and as a result, $tr(G)$ is even.\qed

{\theorem \label{tree-complement} Let $T$ be a tree on $n\geq 4$ vertices which is not a star. Then $$tr(T^c)=\left\lbrace \begin{array}{ll}
	6, & \mbox{ if } T \mbox{ is a bistar.}\\
	5, & \mbox{ if } T \mbox{ is not a bistar.}
	\end{array}\right.$$ }\\
\pf Since $T$ is not a star, $T^c$ is connected. Let $T$ be a tree which is not a bistar. By the arguments used to prove the lower bounds in Theorem \ref{NG-bound}, it follows that $tr(T^c)\geq 5$. Let $d$ and $\overline{d}$ denote the distance in $T$ and $T^c$ respectively. If possible, let $tr(T^c)\geq 6$ and let $u,v,w$ be three vertices attaining the triameter, i.e., $\overline{d}(u,v)+\overline{d}(v,w)+\overline{d}(w,u)=tr(T^c)\geq 6$. If each of $\overline{d}(u,v),\overline{d}(v,w),\overline{d}(w,u)$ is greater than or equal to $2$, then $u,v,w$ forms a triangle in $T$, a contradiction. Thus, at least one of $\overline{d}(u,v),\overline{d}(v,w),\overline{d}(w,u)$ is $1$, say $\overline{d}(u,v)=1$. Then $\overline{d}(v,w)+\overline{d}(w,u)\geq 5$. Thus at least one of $\overline{d}(v,w),\overline{d}(w,u)$ is greater than or equal to $3$, say $\overline{d}(v,w)\geq 3$. 

We claim that $\overline{d}(v,w)= 3$. If possible, let $\overline{d}(v,w)=k\geq 4$. Then there exists vertices $x_1,x_2,\ldots,x_{k-1}$ in $T$ such that $v\sim x_1 \sim x_2 \sim \cdots \sim x_{k-1} \sim w$ is a shortest path joining $v$ and $w$ in $T^c$. Thus, we get a cycle $v\sim x_3 \sim x_1 \sim w \sim x_2 \sim v$ in $T$, a contradiction. Hence $\overline{d}(v,w)= 3$, $\overline{d}(u,w)\geq 2$ and $v\sim x_1 \sim x_2 \sim w$ is a shortest path in $T^c$. 

We claim that $u=x_1$. If possible, let $u\neq x_1$. As $\overline{d}(u,v)=1$, $u\sim v \sim x_1\sim x_2 \sim w$ is a path in $T^c$. Also, $u \not\sim w$ in $T^c$, as $u \sim w$ in $T^c$ implies $\overline{d}(v,w)=2<3$, a contradiction. But this gives rise to a cycle $u \sim w \sim v \sim x_2 \sim u$ in $T$, a contradiction. Thus $x_1=u$ and hence $v\sim u \sim x_2 \sim w$ is a shortest path in $T^c$ and $\overline{d}(u,w)= 2$.

\begin{figure}[h]
	\centering
	\begin{center}
		\includegraphics[scale=.32]{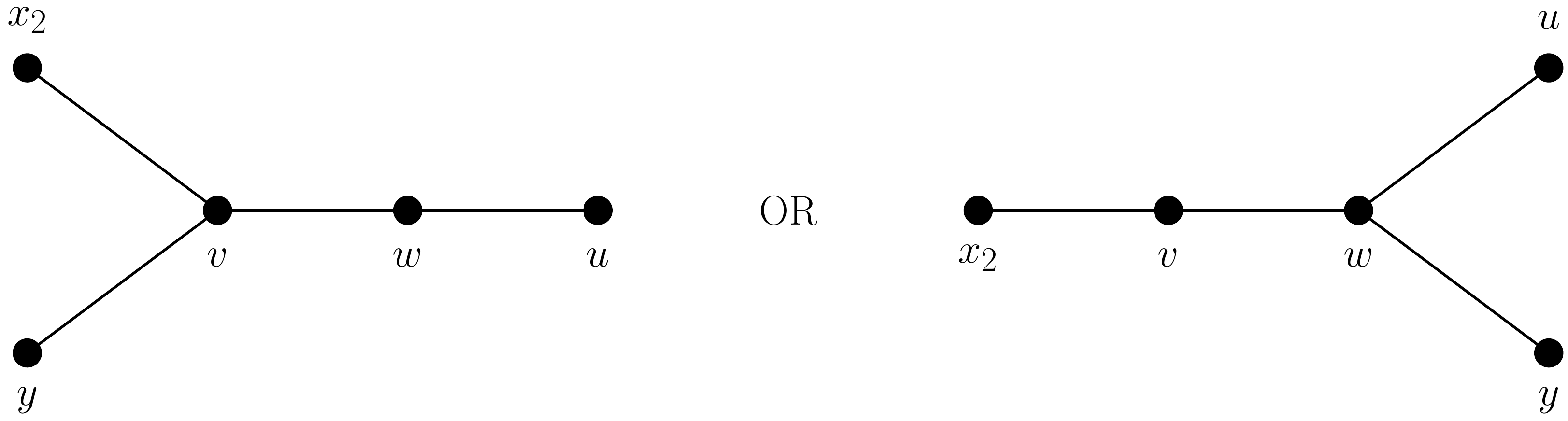}
		\caption{Two Possibilities in Proof of Theorem \ref{tree-complement} }
		\label{bistar-proof-figure-1}	
	\end{center}
\end{figure}

Thus, in $T$, we get a path $x_2\sim v \sim w\sim u$. If $T$ has no other vertex, then $T=P_4$ and it is a bistar with two copies of $K_{1,1}$ joined by an edge. If $T$ has vertices, other than $u,v,w,x_2$, then there exists $y \in V(T)$ such that which is adjacent to exactly one of $u,v,w,x_2$ in $T$ (more than one adjacency creates a cycle in T). If $y$ is adjacent to $x_2$ or $u$ in $T$, we get $v\sim y \sim w$ in $T^c$, i.e, $\overline{d}(v,w)= 2<3$, a contradiction. Thus $y$ is adjacent to $v$ or $w$ in $T$. See Figure \ref{bistar-proof-figure-1}.
If $T$ has no other vertices, $T$ is a bistar with $K_{1,1}$ and  $K_{1,2}$ joined by an edge. If not, then there exists $z \in V(T)$ such that $z$ is adjacent to exactly one of $u,v,w,x_2,y$ in $T$.

We claim that $z$ is adjacent to $v$ or $w$ in $T$. If not, we get a path $v \sim z \sim w$ of length $2$ in $T^c$, i.e., $\overline{d}(v,w)= 2<3$, a contradiction. Thus $z$ is adjacent to $v$ or $w$ in $T$. If $T$ has no other vertices, we again get $T$ to be a bistar. If $T$ has other vertices and as $V(T)$ is finite, in the same way, it can be shown that all other vertices are either adjacent to $v$ or $w$ in $T$. Thus $T$ is a bistar. This is a contradiction to our assumption that $T$ is not a bistar. Hence $tr(T^c)\leq 5$ and hence $tr(T^c)= 5$.

\begin{figure}[h]
	\centering
	\begin{center}
		\includegraphics[scale=.35]{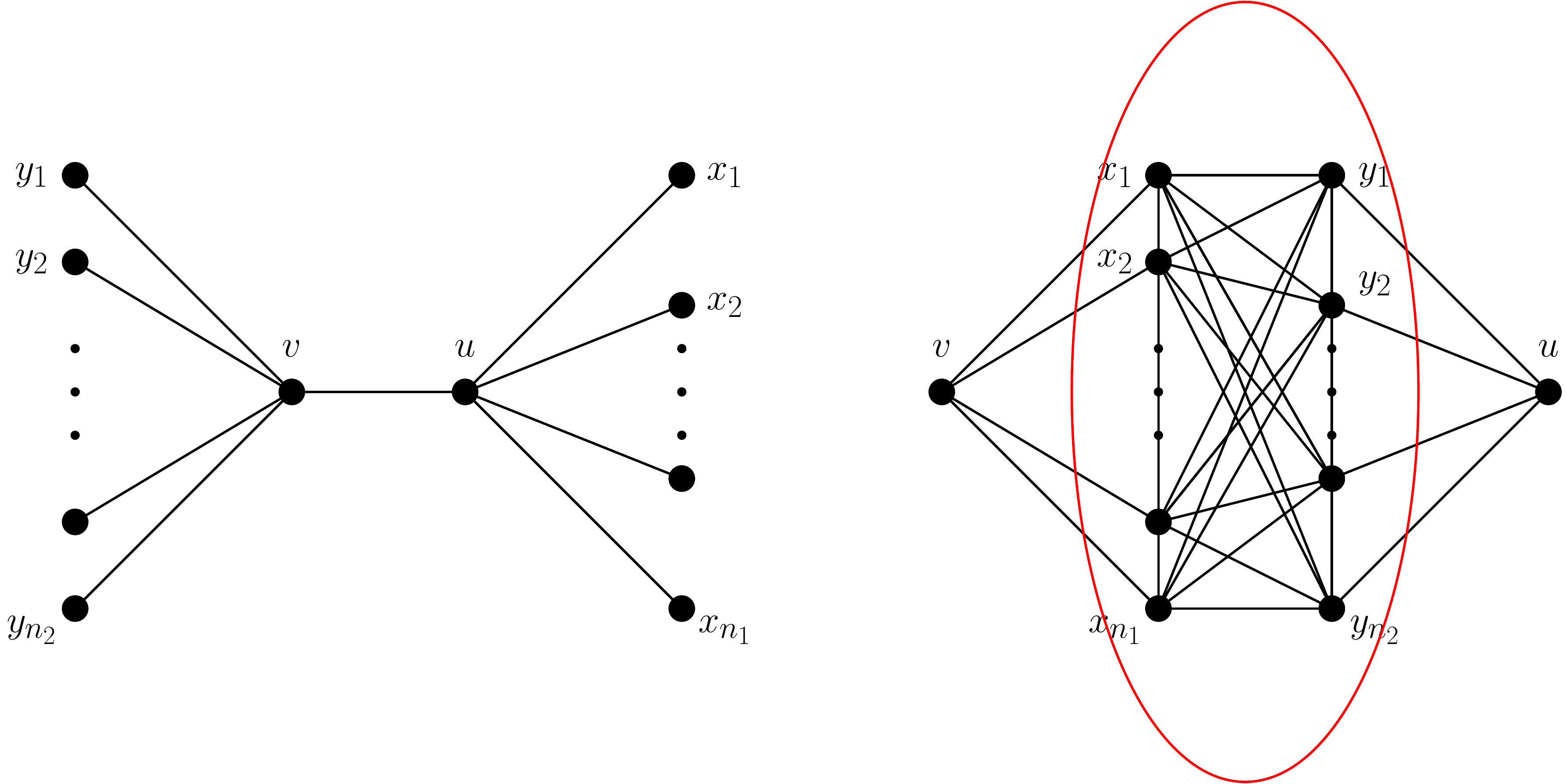}
		\caption{Bistar and its Complement}
		\label{bistar-complement-figure-2}
	\end{center}
\end{figure}

For the other part, let $T=K^{n_1}_{n_2}$ be a bistar as in Figure \ref{bistar-complement-figure-2} (left). Then its complement is as in Figure \ref{bistar-complement-figure-2} (right).

Note that the complement consists of a clique induced by $x_1,x_2,\ldots,x_{n_1},y_1,y_2,\ldots,y_{n_2}$ (indicated in red) and $v$ being adjacent to all the $x_i$'s and $u$ being adjacent to all the $y_i$'s. Thus, for triameter of $T^c$, we need to two of the vertices as $u$ and $v$ and the other to be any one of $x_i$'s or $y_i$'s. Hence, $tr(T^c)=3+1+2=6$.\qed

{\proposition \label{hamiltonian-bound} If $G$ is a Hamiltonian graph on $n$ vertices, then $tr(G)\leq n$.}\\
\pf Since $G$ is a Hamiltonian graph on $n$ vertices, $G$ contains $C_n$ as a subgraph and hence $tr(G)\leq tr(C_n)=n$.\qed

{\theorem \label{vt-bound} If $G=(V,E)$ is a connected vertex transitive graph, then $2\cdot rad(G) \leq tr(G)\leq 3\cdot rad(G)$.}\\
\pf As $G$ is vertex transitive, $V=center(G)=\{x\in V:ecc(x)=rad(G)\}$. Thus, for $u,v,w\in V$, $d(u,v)+d(v,w)+d(w,u)\leq ecc(u)+ecc(v)+ecc(w)=3\cdot rad(G)$. Hence the upper bound follows. The lower bound follows from Corollary \ref{radius-bound}. The tightness of lower and upper bounds follows by taking $G$ as $C_{2n}$ and Petersen graph respectively. \qed

{\theorem \label{srg-bound} If $G$ is a connected strongly regular graph, then $$tr(G)=\left\lbrace \begin{array}{ll}
	5, & \mbox{ if } \overline{G} \mbox{ is triangle-free.}\\
	6, & \mbox{else.}
	\end{array}\right. $$}
\pf Let $G$ be strongly regular with parameters $(n,k,\lambda,\mu)$. Since $G$ is connected,  $\mu\neq 0$ and $G$ is not a complete graph. As a connected strongly regular graph has diameter $2$, $tr(G)\leq 6$. Moreover $\overline{G}$ is again a strongly regular graph with parameter $(n,n-k-1,n-2k+\mu-2,n-2k+\lambda)$. Let $u,v$ be two non-adjacent vertices in $G$, i.e., $d(u,v)=2$. If there exists a vertex $w$ such that $d(u,w)=d(v,w)=2$, then choosing $u,v,w$ as the three vertices we get $tr(G)= 6$. If there does not exist such vertices in $G$, then all vertices other than $u$ and $v$ are either adjacent to $u$ or $v$ or both. Thus, counting the vertices in $G$, we get $n=(k-\mu)+(k-\mu)+\mu +2$, i.e., $n-2k+\mu -2=0$, i.e., $\overline{G}$ is triangle free. In this case, choosing any $w \in N(v)\setminus N(u)$ in $G$, we get $d(u,w)=2$ and $d(v,w)=1$. Then $tr(G)=5$.\qed 

\section{Nordhaus-Gaddum Bounds}
In this section, we prove some Nordhaus-Gaddum type bounds on traimeter of a graph and its complement.

{\lemma Let $G$ be a connected graph such that $\overline{G}$ is connected. Then $tr(G)\geq 7$ implies $tr(\overline{G})\leq 12$. }\\
\pf Since $diam(G)\geq tr(G)/3>2$, it follows that $\gamma(\overline{G})=2$. Thus $tr(\overline{G})\leq 6\gamma(\overline{G})=12$.\qed

{\lemma Let $G=(V,E)$ be a graph such that $G$ and $\overline{G}$ is connected. If $tr(G)>9$, then $tr(\overline{G})\leq 6$.}\\
\pf If possible, let $tr(G)>9$ and $tr(\overline{G})\geq 7$. 
 Let $u,v,w$ be three arbitrary vertices in $V$. 

{\bf Case 1:} If at least one of  $d_{\overline{G}}(u,v),d_{\overline{G}}(v,w),d_{\overline{G}}(w,u)$, say $d_{\overline{G}}(w,u)$ is greater than $1$, then $d_{G}(w,u)=1$. If $d_G(u,v)$ or $d_G(v,w)$ is greater than $3$, then $diam(G)>3\Rightarrow diam(\overline{G})\leq 2 \Rightarrow tr(\overline{G})\leq 6$, a contradiction. Thus $d_G(u,v),d_G(v,w)\leq 3$, i.e., $d_G(u,v)+d_G(v,w)+d_{G}(w,u)\leq 7\leq 9$.

{\bf Case 2:} If  $d_{\overline{G}}(u,v)=d_{\overline{G}}(v,w)=d_{\overline{G}}(w,u)=1$, then $2 \leq d_G(u,v),d_G(v,w),d_G(w,u)\leq 3$ and hence $d_G(u,v)+d_G(v,w)+d_{G}(w,u)\leq 9$.

Combining the two cases we get $tr(G)\leq 9$, which is a contradiction to the assumption and hence the lemma holds. \qed


{\theorem \label{NG-bound} Let $G=(V,E)$ be a graph with $n\geq 4$ vertices such that $G$ and $\overline{G}$ is connected. Then 
	\begin{itemize}
		\item $10\leq tr(G)+tr(\overline{G})\leq 2n+4$,
		\item $25\leq tr(G)\cdot tr(\overline{G})\leq 12(n-1)$ except for a finite family of graphs $\mathcal{F}$,
	\end{itemize} and the bounds are tight.}\\
\pf If $tr(G)>9$, $tr(\overline{G})\leq 6$. Also, by Theorem \ref{order-bound}, $tr(G)\leq 2n-2$ and hence $tr(G)+tr(\overline{G})\leq 2n+4$ and $tr(G)\cdot tr(\overline{G})\leq 12(n-1)$. Let $tr(G)\leq 6$ and if possible, let $tr(G)+tr(\overline{G})> 2n+4$ or $tr(G)\cdot tr(\overline{G})> 12(n-1)$. Then $tr(\overline{G})>2n-2$, a contradiction to Theorem \ref{order-bound}. Thus, if $tr(G)>9$ or $tr(G)\leq 6$, the both the upper bounds hold. Similarly, if $tr(\overline{G})>9$ or $tr(\overline{G})\leq 6$, both the upper bounds hold. 

So the only cases left are when $tr(G),tr(\overline{G})\in \{7,8,9\}$. Thus by Theorem \ref{diameter-bound}, $diam(G)$, $diam(\overline{G})\in \{3,4\}$. However, if $diam(G)$ or $diam(\overline{G})$ equals $4$, then $diam(\overline{G})$ or $diam(G)$ less than or equal to $2$, a contradiction. Thus $diam(G)=diam(\overline{G})=3$.

However, in this cases, for $n\geq 7$, $tr(G)+tr(\overline{G})\leq 18 \leq  2n+4$ and for $n\geq 8$, $tr(G)\cdot tr(\overline{G})\leq 81 \leq 12(n-1)$. 

In \cite{table-six-vertices-connected}, authors provide a complete list of 112 connected graphs on $6$ vertices. Similarly, there are exactly $5$ non-isomorphic graphs (See \cite{list}) on $5$ vertices for which both the graph and its complement is connected. Finally, $P_4$ is the only connected graph on $4$ vertices whose complement is also connected. An exhaustive check (using Sage \cite{sage}) on these graphs revealed that the additive upper bound holds for $n=4,5,6$, and hence the additive upper bound holds for all $n\geq 4$. Also note that for $P_4$, the multiplicative upper bound is an equality.

For the multiplicative upper bound in case of $n=5,6,7$, let us define a family of graphs $\mathcal{F}$ as follows: $$\mathcal{F}=\{G: |G|\in \{5,6,7\}; diam(G)=diam(\overline{G})=3; tr(G),tr(\overline{G})\in \{7,8,9\} \}.$$ From the above discussions, it follows that the multiplicative upper bound holds for all graphs $G$ not in $\mathcal{F}$.

For the lower bounds, observe that as $diam(G)=1$ implies $\overline{G}$ is disconnected, we have $diam(G),diam(\overline{G})\geq 2$, and hence by Theorem \ref{diameter-bound}, $tr(G),tr(\overline{G})\geq 4$. If possible, let $tr(G)=4$, then there exists $u,v,w\in W$, such that $d(u,v)+d(v,w)+d(w,u)=4$. Without loss of generality, let us assume $d(u,v)=2$ and $d(v,w)=d(w,u)=1$. If $G$ is a graph on $3$ vertices, then $P_3$ is the only choice for $G$ satisfying the condition. However, complement of $P_3$ is not connected. Thus we assume that order of $G$ is greater than $3$. Note that for all $z\in V\setminus\{u,v\}$, we have $d(u,z)=d(v,z)=1$ in $G$. But this implies that $\overline{G}$ is disconnected with $u,v$ as one of the components. Thus, to ensure connectedness of $G$ and $\overline{G}$, we have $tr(G),tr(\overline{G})\geq 5$ and hence the additive and multiplicative lower bounds follows.

If $G=P_4$, path on 4 vertices, then $tr(G)=tr(\overline{G})=6$ and hence the upper bounds are tight. If $G=C_5$, cycle on $5$ vertices, then $tr(G)=tr(\overline{G})=5$ and hence the lower bounds are tight. \qed

\begin{figure}[h]
	\centering
	\begin{center}
		\includegraphics[scale=.25]{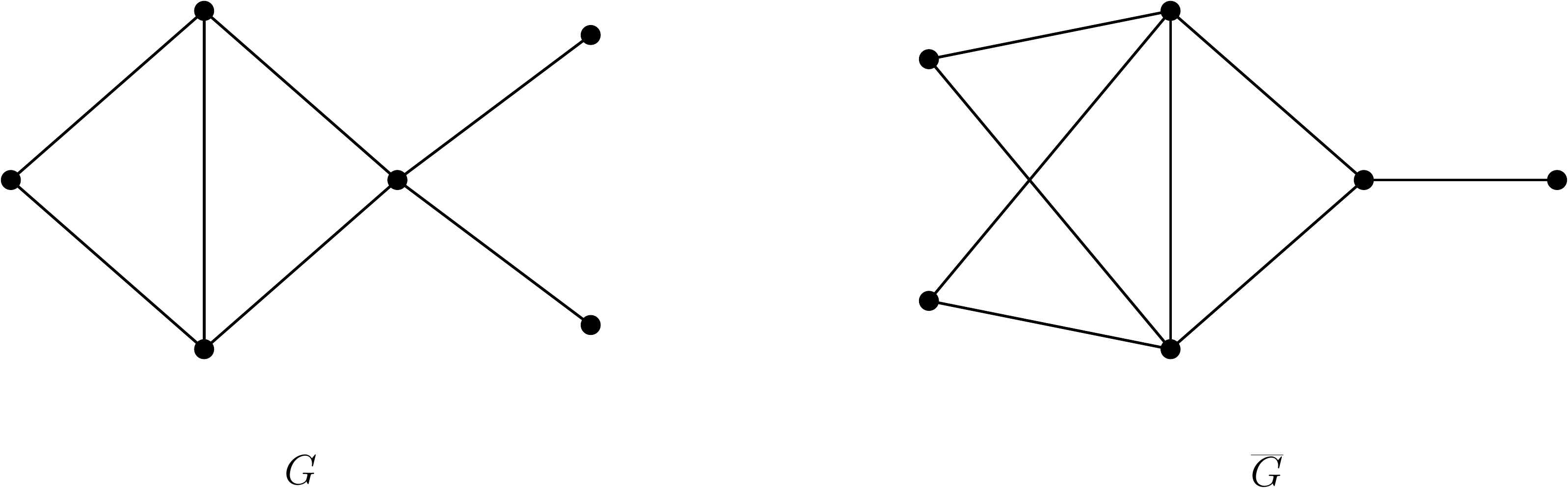}
		\caption{$G,\overline{G}\in \mathcal{F}$}
		\label{F-family-figure}
	\end{center}
\end{figure}

{\remark The multiplicative upper bound may not hold for graphs in $\mathcal{F}$. We demonstrate it in Figure \ref{F-family-figure}. Here $n=6,diam(G)=diam(\overline{G})=3, tr(G)=tr(\overline{G})=8$. Thus $tr(G)\cdot tr(\overline{G})=64>12(6-1)$. }

\section{Conclusion and Open Problems}
In this paper, motivated by a lower bound on radio $k$-coloring in graphs, we formally introduce the idea of triameter in graphs and provide various bounds of various types with respect to other graph parameters. We also provide a shorter proof of a result in \cite{henning-total-lower-bound}. We conclude with two possible directions of further research.

\begin{itemize}
	\item Theorem \ref{(n-1/l-1)-bound} provides a lower bound of $tr(T)$ in terms of its order $n$ and number of leaves $l\geq 3$. Though the bound is tight for $l=3$, the bound loosens as $l$ increases. To find a tighter bound can be an interesting topic of research.
	\item The only lower bound for connected graphs $G$ (not necessarily trees) is in terms of girth (See Theorem \ref{girth-bound}). However, we believe that a better bound is possible in terms of the maximum $\Delta(G)$ and minimum degree $\delta(G)$.
\end{itemize}

\end{document}